\documentclass[12pt,a4paper]{amsart}
\usepackage[utf8]{inputenc}
\usepackage{lmodern}
\usepackage[T1]{fontenc}
\usepackage{graphicx}
\usepackage{amsthm}
\usepackage{amsmath}
\usepackage{amsfonts}
\usepackage{amssymb}
\usepackage[margin=2.50cm]{geometry}
\usepackage[shortlabels]{enumitem}
\usepackage{todonotes}
\usepackage{tikz-cd}
\usepackage{hyperref}
\usepackage{vwcol}
\hypersetup{colorlinks,
	linkcolor={red!50!black},
	citecolor={blue!60!black},
	urlcolor={blue!80!black}
}

\newcommand{\sq}{\subseteq}

\newcommand{\N}{\mathbb{N}}

\newcommand{\R}{\mathbb{R}}

\newcommand{\C}{\mathbb{C}}

\newcommand{\uph}{\upharpoonright}

\newcommand{\A}{\mathcal{A}}

\newcommand{\eps}{\varepsilon}
\newcommand{\M}{\mathcal{M}}
\newcommand{\n}{\mathcal{N}}

\newcommand{\MA}{\mathsf{MA}}
\newcommand{\AD}{\mathsf{AD}}
\newcommand{\PD}{\mathsf{PD}}
\newcommand{\CH}{\mathsf{CH}}
\newcommand{\pp}{\mathcal{P}}
\newcommand{\perpp}{\underline{\perp}}

\newcommand*\diff{\mathop{}\!\mathrm{d}}

\DeclareMathOperator{\ev}{ev}

\DeclareMathOperator{\id}{id}

\DeclareMathOperator{\supp}{supp}

\DeclareMathOperator{\proj}{proj}
\DeclareMathOperator{\extr}{\partial_e}
\DeclareMathOperator{\tr}{tr}
\DeclareMathOperator{\aut}{Aut}

\theoremstyle{plain}
\newtheorem{theorem}{Theorem}[section]
\newtheorem{lemma}[theorem]{Lemma}
\newtheorem{proposition}[theorem]{Proposition}
\newtheorem{corollary}[theorem]{Corollary}
\newtheorem{question}[theorem]{Question}

\newtheorem{claim}[theorem]{Claim}
\newtheorem{fact}[theorem]{Fact}

\theoremstyle{definition}

\newtheorem{remark}[theorem]{Remark}
\newtheorem{example}[theorem]{Example}

\title{Orthogonality of measures and states}
\author{Severin Mejak}
\address{Department of Mathematical Sciences, University of Copenhagen, Universitetsparken 5, 2100 Copenhagen, Denmark}
\email{seme@math.ku.dk}

\begin{document}

    \maketitle
    
    \begin{abstract}
    We give a short proof of the theorem due to Preiss and Rataj stating that there are no analytic maximal orthogonal families (mofs) of Borel probability measures on a Polish space. When the underlying space is compact and perfect, we show that the set of witnesses to non-maximality is comeagre. Our argument is based on the original proof by Preiss and Rataj, but with significant simplifications. The proof generalises to show that under $\mathsf{MA} + \neg \mathsf{CH}$ there are no $\mathbf{\Sigma^1_2}$ mofs, that under $\mathsf{PD}$ there are no projective mofs and that under $\mathsf{AD}$ there are no mofs at all. We also generalise a result due to Kechris and Sofronidis, stating that for every analytic orthogonal family of Borel probability measures there is a product measure orthogonal to all measures in the family, to states on a certain class of C*-algebras.
    \end{abstract}
    
    \section*{Introduction}
    In this paper we consider orthogonality, first for Borel probability measures on Polish spaces and then for states on separable C*-algebras. In 1985 Preiss and Rataj proved the following theorem with $X = [0, 1]$, see \cite{preiss_rataj}.
    
    \begin{theorem}\label{thm_measures}
        Suppose that $X$ is an uncountable Polish space. Then there is no analytic maximal orthogonal family of Borel probability measures on $X$.
    \end{theorem}
    
    This answered an open question from \cite{mauldin}. The idea of the proof from \cite{preiss_rataj} is to use a Baire category argument. However, once the authors prepared the scene for the application of the Baire category theorem, they resorted to a couple of technical lemmas, which relied on restricting Borel probability measures on $[0, 1]$ to finite unions of closed subintervals. For the proof of one of the lemmas they also employed Banach--Mazur games. Consequently,
    the question whether there is a shorter and simpler proof remained open.
    
    In 1999, Kechris and Sofronidis (see Thoerem 3.1 in \cite{strongregularity}) found an alternative short proof which uses the theory of \emph{turbulence} (see \cite{hjorth} for a great introduction to turbulence). As part of their proof, they defined an embedding of the Cantor space $2^\N$ into the space of Borel probability measures (using the work of Kakutani from \cite{kakutani}), assigning to every $x \in 2^\N$ a product measure $\mu_{\alpha(x)}$. They proved that for every analytic orthogonal family, there is some $x \in 2^\N$ so that $\mu_{\alpha(x)}$ is a witness to non-maximality. Their proof has as a consequence that the relation $\sim$ of measure equivalence between Borel probability measures is not classifiable by countable structures.
    
    Almost two decades later Schrittesser and T\"{o}rnquist used the same embedding of $2^\N$ into the space of measures to prove (see Theorem 5.5 of \cite{mrl}) that an argument using a weaker form of turbulence suffices to prove Theorem \ref{thm_measures}. Since the theory of turbulence requires some background knowledge, one might argue that even thought the proofs from \cite{strongregularity} and \cite{mrl} are \emph{shorter}, they are not necessarily \emph{simpler}.
    
    In this article, we first go back to the original idea of Preiss and Rataj to use a Baire category argument to prove Theorem \ref{thm_measures}. We were able to use the Kuratowski--Ulam theorem and some elementary convexity theory, to give a short and straightforward proof of Theorem \ref{thm_measures}. The argument works to show the following strengthening, where for $\A \sq P(X)$ (here $P(X)$ denotes the space of Borel probability measures on $X$), we let $\A^\perp := \{\nu \in P(X) : (\forall \mu \in \A)\, \nu \perp \mu\}$.
    
    \begin{theorem}\label{thm_measures_comeagre}
        Suppose that $X$ is a compact perfect Polish space. Then for every analytic orthogonal family $\A \sq P(X)$, the set $\A^\perp$ is comeagre. In particular, when $\A \sq P(2^\N)$ is a $\Sigma^1_1$ orthogonal family, there is a $\Delta^1_1$-witness to non-maximality.
    \end{theorem}
    
    Actually, under additional assumptions, our method yields the following.
    
    \begin{theorem}\label{MA_AD_measures}
    Suppose that $X$ is an uncountable Polish space.
            \begin{enumerate}
                \item Assume $\MA$ and $\neg \CH$. Then no $\mathbf{\Sigma^1_2}$ orthogonal family $\A \sq P(X)$ is maximal.
                \item Assume $\PD$. Then no projective orthogonal family $\A \sq P(X)$ is maximal.
                \item Assume $\AD$. Then no orthogonal family $\A \sq P(X)$ is maximal.
            \end{enumerate}
    If moreover $X$ is compact perfect, then in each of the above cases $\A^\perp$ is comeagre.
    \end{theorem}
    
    It is well-known that via the Riesz--Markov--Kakutani representation theorem, Borel probability measures on a compact Polish space $X$ are precisely states on the commutative C*-algebra of complex-valued continuous functions on $X$. In \cite{dye}, Dye introduced the notion of absolute continuity for states on C*-algebras. Being a pre-order, it naturally gives rise to a notion of orthogonality. However, it turns out that this notion is ill-behaved even for states on the matrix algebra $M_2(\C)$.
    
    There is another natural notion of orthogonality for states, which we call strong orthogonality and denote by $\perpp$. This notion of orthogonality shares many nice properties with orthogonality of measures, with which it coincides when the C*-algebra is commutative. Hence it is natural to ask ourselves whether Theorem \ref{thm_measures} holds for non-commutative separable unital C*-algebras and strong orthogonality as well.
    
    Since the original proof by Preiss and Rataj relied on restrictions of measures to compact subspaces, it is not clear how to generalise that proof. The idea from the proof of Theorem \ref{thm_measures} seems more promising, but there are still some steps for which we do not know if they hold for strong orthogonality for states.
    
    On the other hand, it turns out that the idea of Kechris and Sofronidis from \cite{strongregularity} can easily be extended to a class of separable unital C*-algebras.
    
    \begin{theorem}\label{thm_product_states}
        Suppose $A$ is a separable unital C*-algebra, which contains a copy of $C(2^\N)$ as a subalgebra and for which there is a conditional expectation $E: A \to C(2^\N)$. Then for every strongly orthogonal $\A \sq S(A)$ there is $\alpha \in (0, 1)^\N$ so that $\Tilde{\mu}_\alpha \perpp \psi$ for every $\psi \in \A$, where $\Tilde{\mu}_\alpha$ is the extension of the state, corresponding to the product measure
        \[
        \prod_{n \in \N} (\alpha(n) \delta_0 + (1-\alpha(n)) \delta_1),
        \]
        from $C(2^\N)$ to $A$.
    \end{theorem}
    
    As in \cite{strongregularity}, along the way of proving this theorem we also get that for C*-algebras $A$, satisfying the assumptions of the theorem, the relation $\sim$ on $S(A)$ is not classifiable by countable structures.
    
    Natural examples of C*-algebras, for which the assumptions of Theorem \ref{thm_product_states} are satisfied, include the CAR algebra $M_{2^\infty}$ and the Cuntz algebra $\mathcal{O}_2$. Moreover, for any $A$ satisfying assumptions of Theorem \ref{thm_product_states} also the reduced crossed product $A \rtimes_{\alpha, r} \Gamma$ (for any countable discrete group $\Gamma$ and any homomorphism $\alpha: \Gamma \to \aut(A)$) and the tensor product $A \otimes B$ (for any separable unital C*-algebra $B$) satisfy the assumptions of Theorem \ref{thm_product_states}.
    
     In 1969 Bures (see \cite{bures}) proved an extension of Kakutani's result from \cite{kakutani} to semi-finite von Neumann algebras. Instead of absolute continuity and orthogonality between states, Bures considered when two product states give rise to isomorphic tensor products of von Neumann algebras. This was extended to all von Neumann algebras by Promislow in \cite{promislow}.
     
     As a consequence of the main ingredient of the proof of Theorem \ref{thm_product_states}, we get the following version of Kakutani's theorem for states, involving absolute continuity and strong orthogonality.
     
     \begin{proposition}\label{kakutani_states}
        Suppose that $(\alpha_n)_{n \in \N}, (\beta_n)_{n \in \N} \in [\frac{1}{4}, \frac{3}{4}]^\N$ and let
    \[\phi_n := \alpha_n \ev_{1, 1} + (1-\alpha_n) \ev_{2, 2} \quad \text{and} \quad  \psi_n := \beta_n \ev_{1, 1} + (1 -\beta_n) \ev_{2, 2}
    \]
    be states on $M_2(\C)$. Let also $\phi := \otimes_{n = 0}^\infty \phi_n$ and $\psi := \otimes_{n = 0}^\infty \psi_n$ be the product states on $M_{2^\infty}$. Then in $S(M_{2^\infty})$, either $\phi \sim \psi$ or $\phi \perpp \psi$ according to whether
        \[
        \sum_{n \in \N} (\alpha_n - \beta_n)^2
        \]
    converges or diverges respectively.
     \end{proposition}
    
    \subsection*{Structure of the paper}
    
    The paper aims to target interested readers from descriptive set theory, measure theory and C*-algebras. Due to different backgrounds, we try to give as much details as possible and add references to literature containing more information about the discussed topics. Readers not familiar with set-theoretic notions such us $\MA$ or $\mathbf{\Sigma^1_2}$, can skip the parts where we consider them, with no effect to understanding the rest of the paper.
    
    In section \ref{measures}, we give proofs of Theorems \ref{thm_measures}, \ref{thm_measures_comeagre} and \ref{MA_AD_measures}. This is followed by section \ref{states}, where we first present absolute continuity and two notions of orthogonality for states. Subsection \ref{product_states} recalls the idea and some notions from \cite{strongregularity} and proves Theorem \ref{thm_product_states} and Proposition \ref{kakutani_states}. We conclude the paper by discussing related topics and listing some open problems in section \ref{open_problems}.
    
    \subsection*{Acknowledgments}
    
    The author is above all deeply grateful to his PhD advisor Asger T\"{o}rnquist for expressing his belief that the proof of Theorem \ref{thm_measures} from \cite{preiss_rataj} can probably be simplified and that the idea of \cite{strongregularity} can likely be used to prove a similar statement for states. I am also thankful to Asger T\"{o}rnquist for numerous discussions on the topic which contributed heavily to this paper.
    
    The author is indebted to Ryszard Nest for discussions on generalisations of notions from measures to states on C*-algebras, without which the general C*-algebra part of the article would not be possible. I am thankful to Magdalena Musat and Mikael R{\o}rdam for pointing out the classes of C*-algebras satisfying the assumptions of Theorem \ref{thm_product_states}. The author also thanks Alexander Frei for valuable conversations about C*-algebras and Forte Shinko for discussions about measures.
    
    During later stages of writing this paper, the author was visiting Caltech and is greatly appreciative to Alexander Kechris and other members of the Caltech logic group for hosting him. I am also thankful to Alexander Kechris for comments on a draft version of this paper, and for suggesting questions appearing in subsection \ref{nice_bad}.\\
    
    \begin{center}
    \begin{vwcol}[widths={0.18, 0.82}, justify= flush, rule = 0pt]
    \hspace{0.5cm}
	\includegraphics{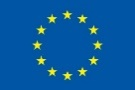}
	\vfill\eject
	\phantom{ }\hspace{-0.6cm}
	\vspace{-0.3cm}
	\phantom{ }\\
	\scriptsize{\noindent This project has received funding from the European Union’s Horizon 2020 research and\\innovation programme under the Marie Sk\l{}odowska-Curie grant agreement No 801199}
    \end{vwcol}
    \end{center}
    
    \section{Borel probability measures}\label{measures}
    
    In this section we give a very short proof of the classical result due to Preiss and Rataj by simplifying some steps of their proof from \cite{preiss_rataj}. The setup of using a Baire category argument is the same, what is new is that we replace the technical part of the proof from \cite{preiss_rataj}, which uses extensions of measures defined on subspaces with a more straightforward argument. We start by recalling some basic properties of Borel probability measures.
    
    Let $X$ be a Polish space. We denote by $C(X)$ the set of continuous complex-valued functions on $X$. With $P(X)$, we denote the collection of Borel probability measures on $X$, endowed with the topology generated by maps $\mu \mapsto \int f d\mu$ for $f$ ranging over $C_b(X, \R) := \{f: X \to \R \, | \, f \text{ continuous and bounded}\}$. Recall that when $X$ is Polish, then so is $P(X)$, and that if moreover $X$ is compact, so is $P(X)$. See section 17.E of \cite{kechris} for more about $P(X)$.
    
    For two Borel measures $\mu, \nu$ on $X$, we denote $\mu$ being \emph{absolutely continuous} with respect to $\nu$ (i.e., for every Borel subset $B \sq X$, if $\nu(B) = 0$, then $\mu(B) = 0$) by $\mu \ll \nu$. We say that $\mu, \nu \in P(X)$ are \emph{measure equivalent}, denoted by $\mu \sim \nu$, if $\mu \ll \nu \land \nu \ll \mu$ and that $\mu, \nu \in P(X)$ are \emph{orthogonal} (another term often used is \emph{singular}), denoted by $\mu \perp \nu$, if there is no $\rho \in P(X)$ with $\rho \ll \mu$ and $\rho \ll \nu$. Observe that $\mu \perp \nu$ is equivalent to existence of a Borel $B \sq X$ with $\mu(B) = 1$ and $\nu(B) = 0$. Recall (see e.g. \cite{strongregularity}) that $\ll$ satisfies the \emph{ccc-below property}, i.e., for every $\mu \in P(X)$ there is no uncountable family $\{\nu_i : i \in I\} \sq P(X)$ with the property that for $i \neq j \in I$ it holds that $\nu_i \perp \nu_j$ and $\nu_i \ll \mu$.
    
    For a signed Borel measure $\sigma$ on $X$, 
    \[
    ||\sigma|| := \sup \{ |\sigma(B)| : B \sq X \text{ Borel}\}
    \]
    defines a norm on the space of signed Borel measures. Then the map $P(X) \times P(X) \to \R_{\geq 0}$, defined by $(\mu, \mu) \mapsto ||\mu - \nu||$ is lower semicontinuous, when $P(X)$ is equipped with the Polish topology defined above. Consequently, for $\eps \in (0, 1)$ the set
    \[
    \{(\mu, \nu) \in P(X) \times P(X) : ||\mu - \nu|| < \eps\}
    \]
    is $F_\sigma$ (i.e., it is a countable union of closed sets). It is also immediate to see that for $\mu, \nu \in P(X)$ we have
    \[
    \mu \perp \nu \text{ if and only if } ||\mu - \nu|| = 1,
    \]
    so that the relation $\mu \perp \nu$ is $G_\delta$ (i.e., a countable intersection of open sets). We continue with a lemma from convexity theory.
    
    \begin{lemma}\label{convex_comb}
        Suppose that $V$ is an open convex subset of a locally convex topological vector space $E$. Then the map $V \times V \times [0, 1] \to V$, defined by $(x, y, t) \mapsto t x + (1-t) y$ is continuous and open.
    \end{lemma}
    
    \begin{proof} 
        Continuity holds because $E$ is a topological vector space. To check that the map is open, take $U_0, U_1 \sq V$ convex open, $O \sq [0, 1]$ convex open and $x, y, t$ in $U_0, U_1, O$ respectively. Then let
        \[U_2 := ((1-t) (y - x) + U_0) \cap (t (x - y) + U_1),\]
        which is convex open and contains $t x + (1-t) y$. Clearly $U_2$ is contained in the image of $U_0 \times U_1 \times O$, hence this completes the proof that the map is open.
    \end{proof}
    
    We are now ready to prove Theorem \ref{thm_measures}, which we restate for reader's convenience.
    
    \begingroup
    \def\thetheorem{\ref{thm_measures}}
    \begin{theorem}
        Suppose that $X$ is an uncountable Polish space. Then there is no analytic maximal orthogonal family of Borel probability measures on $X$.
    \end{theorem}
    \addtocounter{theorem}{-1}
    \endgroup

    \begin{proof}
        Suppose for contradiction that $\A$ is an analytic maximal pairwise orthogonal family of Borel probability measures on $X$. Observe that we can without loss of generality assume that $X$ is perfect and compact. Indeed, for a general uncountable Polish space $X$, there is a compact perfect subspace $Y$ of $X$ (see Theorems 6.4 and 6.2 of \cite{kechris}). Then
        \[
        \A' := \left\{\frac{1}{\mu(Y)} \mu \uph Y : \mu \in \A \land \mu(Y) > 0  \right\}
        \]
        is clearly a maximal analytic family of pairwise orthogonal Borel probability measures on the uncountable perfect compact Polish space $Y$.
        
        The first steps follow the proof of Theorem \ref{thm_measures} from \cite{preiss_rataj}. We include these steps with a little more detail for reader's convenience. For every $k \in \N$ denote by $E_k$ the space of $k$-element subsets of $\A$, equipped with the usual topology in which it is clearly analytic. Fix some $\eps \in (0, 1)$ and define
        \[
        H_{k, \eps} := \{\nu \in P(X) : (\exists F \in E_k) \, (\forall \mu \in F)\, ||\nu - \mu || < \eps\},
        \]
        which is evidently analytic, and for a fixed $\tau \in (0, \eps)$ define also
        \[
        U_{k, \eps}^\tau := H_{k, \eps - \tau} \setminus H_{k + 1, \eps},
        \]
        which thus has the property of Baire. Since $\A$ is maximal orthogonal, we have that for every $\nu \in P(X)$ there is some $\mu \in \A$ with $||\nu - \mu|| < 1$. Moreover, since $\A$ consists of pairwise orthogonal measures, it holds that for any $0 \leq \sigma < 1$ and any $\nu \in P(X)$ the set
        \[
        \{\mu \in \A : ||\nu - \mu|| < \sigma\}
        \]
        is finite. Indeed, if it were infinite find $n \geq 1$, such that $1 - \sigma > 1/n$. Then there are some $\mu_0, \ldots, \mu_n \in \A$ with $||\nu - \mu_j|| < \sigma < 1 - 1/n$ for all $0 \leq j \leq n$. Therefore there are pairwise disjoint Borel subsets $D_0, \ldots, D_n$ of $X$ such that for all $0 \leq i, j \leq n$, $i \neq j$ it holds that $\mu_i(D_i) = 1$ and $\mu_i(D_j) = 0$. Thus we have for all $0 \leq i \leq n$ that $\nu(D_i) > 1/n$. But then
        \[
        \nu(X) \geq \nu\left(\bigcup_{i = 0}^n D_i\right) > (n + 1) \frac{1}{n} > 1,
        \]
        which is of course a contradiction.
        
        It is clear then that every $\nu \in P(X)$ is in some $U_{k, 1/n}^{1/m}$ for some $k \geq 1, n > 1$ and $m > n$. Hence we have that
        \[
        P(X) = \bigcup_{k \geq 1} \bigcup_{n > 1} \bigcup_{m > n} U_{k, 1/n}^{1/m},
        \]
        and since $P(X)$ is a Baire space, it must hold that for some $k, \eps := 1/n$ and $\tau := 1/m$ it holds that $U_{k, \eps}^\tau$ is comeagre in a non-empty convex open set $V \sq P(X)$ (we can assume convexity of $V$, since $P(X)$ is locally convex). From now on, our proof diverges from the path taken in \cite{preiss_rataj}.
        
        \begin{claim}\label{Kuratowski-Ulam}
        There is $\nu \in U^\tau_{k, \eps} \cap V$ and $C \sq U_{k, \eps}^\tau \cap V$, which is comeagre in $V$, so that for every $\mu \in C$ the set
        \[
        M_\mu := \{t \in [0, 1] : t \nu + (1-t) \mu \in U^\tau_{k, \eps} \}
        \]
        is comeagre in $[0, 1]$.
        \end{claim}
        
        \noindent
        \emph{Proof of Claim.} The map $V \times V \times [0, 1] \to V$, defined by $(\nu, \mu, t) \mapsto t \nu + (1- t) \mu$ is continuous and open by Lemma \ref{convex_comb}, so the preimage of $U_{k, \eps}^{\tau} \cap V$ under this map is also a comeagre subset of $V \times V \times [0, 1]$. Now the Kuratowski-Ulam theorem (see iii of Theorem 8.41 from \cite{kechris}) implies the desired result. \hfill $\dashv$
        
        $ $ \newline
        Let $\nu$ and $C$ be as in the claim above. Next, we introduce the following notation: for $\rho \in U_{k, \eps}^{\tau}$ write
        \[
        N_\rho := \{\mu \in \A : ||\rho - \mu || < \eps - \tau\},
        \]
        and by definition of $U_{k, \eps}^\tau$ observe that $N_\rho = \{\mu \in \A : ||\rho - \mu || < \eps\}$ and has precisely $k$ elements.
        
        \begin{claim}\label{constant}
            For every $\rho \in C$ we have that $N_\nu = N_\rho$.
        \end{claim}
        
        \noindent
        \emph{Proof of Claim.} Fix any $\rho \in C$ and put
        \[
        T:= \{t \in M_\rho : N_{t \nu + (1-t) \rho} = N_\nu\}.
        \]
        and $s := \sup T$. We will show that $s = 1$ and that $s \in T$.
        
        First, suppose for contradiction that $s < 1$. Find $t \in T$, such that $s - t < \tau/5$. Since $M_\rho \sq [0, 1]$ is comeagre, find $u \in M$ with $u \geq s$ and $u - s < \tau/5$. Then 
        \begin{align*}
            ||(u \nu + (1-u) \rho) - (t \nu + (1-t) \rho)|| &= || (u - t) \nu + (t - u) \rho||\\
            &\leq |u-t|\, ||\nu|| + |u-t|\, ||\rho|| = 2 (u-t) < \frac{4 \tau}{5}.
        \end{align*}
        This implies that $N_{u \nu + (1-u) \rho} = N_{t \nu + (1-t) \rho} = N_\nu$, and thus $u \in T$, which is a contradiction.
        
        So we have that $s = 1$. Find $t \in T$ with $|1 - t| < \tau/3$, and observe by the same reasoning as before that since $1 \in M_\rho$, we have that $N_\rho = N_{t \nu + (1-t) \rho} = N_\nu$, completing the proof. \hfill $\dashv$
        
        $ $ \newline
        Let now $\mu_0, \ldots, \mu_{k-1} \in \A$ be such that $N_\nu = \{\mu_0, \ldots, \mu_{k-1}\}$. But since it holds that for any $\mu \in P(X)$ the set $\mu^\perp := \{\rho \in P(X) : \mu \perp \rho \} $ is comeagre (see Proposition 4.1 from \cite{strongregularity} and note that this is where we need that $X$ is perfect compact), also the set
        \[
        B := \bigcap_{j = 0}^{k-1}\, \mu_j^\perp
        \]
        is comeagre and in particular comeagre in $V$. So both $B$ and $C$ are comeagre in the open set $V$, which is of course a contradiction.
        \end{proof}
        
        For a (pairwise orthogonal) family $\A \sq P(X)$, observe that the set of witnesses to non-maximality
        \[
        \A^\perp := \{\nu \in P(X) : (\forall \mu \in \A) \, \nu \perp \mu\}
        \]
        is co-analytic (in particular, it has the Baire property) and by Theorem \ref{thm_measures} it is non-empty. When $X$ is a perfect compact Polish space, we have the following strengthening.
        
    \begingroup
    \def\thetheorem{\ref{thm_measures_comeagre}}
    \begin{theorem}
            Suppose that $X$ is a compact perfect Polish space. Then for every analytic orthogonal family $\A \sq P(X)$, the set $\A^\perp$ is comeagre. In particular, when $\A \sq P(2^\N)$ is a $\Sigma^1_1$ orthogonal family, there is a $\Delta^1_1$-witness to non-maximality.
        \end{theorem}
    \addtocounter{theorem}{-1}
    \endgroup
        
        \begin{proof}
            Suppose for contradiction that $\A^\perp$ is not comeagre. Then there is a non-empty convex open set $O \sq P(X)$, in which $\A^\perp$ is meagre. Let $Z$ be a dense (in $O$) $G_\delta$ subset of $O \setminus \A^\perp$; in particular, $Z$ is a Polish subspace of $P(X)$. Note that $\A$ is maximal orthogonal in $Z$, i.e., for every $\nu \in Z$ there is some $\mu \in \A$ with $||\nu - \mu || < 1$. Now we follow the proof of Theorem \ref{thm_measures}, but this time we use the Baire category theorem in $Z$.
            
            We use the above defined $E_k$ and redefine $H_{k, \eps}$ and $U^\tau_{k, \eps}$ as follows. Fix some $\eps \in (0, 1)$ and set
            \[
            H_{k, \eps} := \{\nu \in Z : (\exists F \in E_k) \, (\forall \mu \in F)\, ||\nu - \mu || < \eps\},
            \]
            which is clearly analytic, and for a fixed $\tau \in (0, \eps)$ define also
            \[
            U_{k, \eps}^\tau := H_{k, \eps - \tau} \setminus H_{k + 1, \eps},
            \]
            which thus has the property of Baire. It is clear then that every $\nu \in Z$ is in some $U_{k, 1/n}^{1/m}$ for some $k \geq 1, n > 1$ and $m > 1$. Hence we have that
            \[
            Z = \bigcup_{k \geq 1} \bigcup_{n > 1} \bigcup_{m > n} U_{k, 1/n}^{1/m},
            \]
            and so it must hold that for some $k, \eps := 1/n$ and $\tau := 1/m$ the set $U_{k, \eps}^\tau$ is comeagre in a non-empty convex open set $V \sq O$. The proofs of Claims \ref{Kuratowski-Ulam} and \ref{constant} work without changes to show the following two claims respectively.
        
            \begin{claim}
            There is $\nu \in U^\tau_{k, \eps} \cap V$ and $C \sq U_{k, \eps}^\tau \cap V$, which is comeagre in $V$, so that for every $\mu \in C$ the set
            \[
            M_\mu := \{t \in [0, 1] : t \nu + (1-t) \mu \in U^\tau_{k, \eps} \}
            \]
            is comeagre in $[0, 1]$.
            \end{claim}
        
            \begin{claim}
            For every $\rho \in C$ we have that $N_\nu = N_\rho$.
            \end{claim}
        
            \noindent
            Finally, let $\mu_0, \ldots, \mu_{k-1} \in \A$ be such that $N_\nu = \{\mu_0, \ldots, \mu_{k-1}\}$. Again, the set
            \[
            B := \bigcap_{j = 0}^{k-1}\, \mu_j^\perp
            \]
            is comeagre and in particular comeagre in $V$. So both $B$ and $C$ are comeagre in the open set $V$, which is again a contradiction.
            
            For the ``in particular'' part of the theorem, note first that $P(2^\N)$ is a recursively presentable Polish space (see \cite{moschovakis} for the definition of the notion and \cite{fischer_tornquist} for why $P(2^\N)$ is recursively presentable), so it makes sense to talk about lightface pointclasses in $P(2^\N)$. To get a $\Delta^1_1$-witness to non-maximality of a $\Sigma^1_1$ orthogonal family $\A \sq P(2^\N)$, use Corollary 4.1.2 of \cite{effective_kechris} on $\A^\perp$, which, as we have just proved, is a comeagre $\Pi^1_1$ set.
        \end{proof}
        
        In \cite{martin_solovay}, Martin and Solovay show that if Martin's axiom ($\MA$) holds and Continuum hypothesis ($\CH$) fails, then all $\mathbf{\Sigma^1_2}$ sets of reals have the Baire property. Recall also that the Axiom of projective determinacy ($\PD$) implies that all projective sets of reals have the Baire property and that the Axiom of determinacy ($\AD$) implies that all sets of reals have the Baire property (see e.g. Theorem 33.3 in \cite{jech}). It is clear that we can substitute sets of reals with subsets of $P(X)$ for a Polish space $X$.
        
    \begingroup
    \def\thetheorem{\ref{MA_AD_measures}}
    \begin{theorem}
    Suppose that $X$ is an uncountable Polish space.
            \begin{enumerate}
                \item Assume $\MA$ and $\neg \CH$. Then no $\mathbf{\Sigma^1_2}$ orthogonal family $\A \sq P(X)$ is maximal.
                \item Assume $\PD$. Then no projective orthogonal family $\A \sq P(X)$ is maximal.
                \item Assume $\AD$. Then no orthogonal family $\A \sq P(X)$ is maximal.
            \end{enumerate}
    If moreover $X$ is compact perfect, then in each of the above cases $\A^\perp$ is comeagre.
    \end{theorem}
    \addtocounter{theorem}{-1}
    \endgroup
        
    \begin{proof}
        Repeat the proof of Theorem \ref{thm_measures} (or  \ref{thm_measures_comeagre} in case $X$ is compact perfect), using the respective assumed axiom to get that the sets $U^\tau_{k, \eps}$ have the Baire property. The rest of the proof is the same.
    \end{proof}
        
    \begin{remark}
        Let $a \in \N^\N$ and note that $\omega_1^{L[a]} < \omega_1$ implies that for every $\Sigma^1_2[a]$ orthogonal family $\A \sq P(2^\N)$, the set $\A^\perp$ is comeagre. The reason is again that $\omega_1^{L[a]} < \omega_1$ implies that every $\Sigma^1_2[a]$ subset of $P(2^\N)$ has the Baire property (see Corollary 14.3 from \cite{kanamori}).
    \end{remark}    
        
    \section{States on separable C*-algebras}\label{states}
    
    \subsection{Absolute continuity and orthogonality for states}
    
    All C*-algebras considered will be unital and separable. For a C*-algebra $A$, we denote by $S(A)$ the collection of states on $A$ (i.e., positive linear functionals, which map the unit of $A$ to 1) and by $PS(A)$ the collection of pure states on $A$ (i.e., the states that are the extreme points of the compact convex set $S(A)$). With $\proj(A)$ we denote the set of projections in $A$. See e.g. \cite{blackadar} and \cite{brown_ozawa} for other standard notions from C*-algebras theory.
    
    It is well-known that every commutative C*-algebra $A$ is *-isomorphic (via the Gelfand transform) to $C(M_A)$, where $M_A$ is the maximal ideal space of $A$ (which can in turn be described as the space of characters on $A$ (i.e., non-zero algebra homomorphisms from $A$ to $\C$)). Furthermore, $M_A$ is compact Polish, being contained in $B_1(A^*)$. So when considering commutative separable C*-algebras we can restrict our attention to $C(X)$ for $X$ compact Polish.
    
    By Riesz--Markov--Kakutani representation theorem we know that $S(C(X))$ is actually the same as $P(X)$, and indeed, this is how $P(X)$ got its topology. Note here that a state $\phi \in S(C(X))$ is determined by its values on real-valued functions (real-valued functions are the self-adjoint elements in $C(X)$).
    
    So it is natural to try to generalise notions from measure theory to states on C*-algebras. In \cite{dye}, Dye defined the notion of absolute continuity for states on $\sigma$-finite von Neumann algebras and proved a version of Radon-Nikodym theorem. There is an abundance of alternative formulations of absolute continuity for states on C*-algebras, some of them equivalent to the one given here, some weaker and some stronger. Standard results from measure theory like Lebesgue decomposition theorem generalise to states (this is the content of \cite{dye}). See also e.g. \cite{henle_1972}, \cite{hiai} or \cite{inoue} for some different formulations of absolute continuity for states and various results on generalisations.
    
    For any \emph{pre-order}  (i.e., a reflexive and transitive relation) $\preccurlyeq$ on a set $X$ one says that $x, y \in X$ are \emph{orthogonal}, in symbols $x \perp y$, if there is no $z \in X$ with $z \preccurlyeq x, y$. Note that orthogonality of measures from the previous section is just orthogonality associated with the pre-order $\ll$. Accordingly, we define orthogonality for states, associated to Dye's notion of absolute continuity.
    
    It is natural to ask whether Theorem \ref{thm_measures} can be generalised to non-commutative unital C*-algebras $A$. As we will see, the orthogonality relation associated with Dye's notion of absolute continuity in general does not satisfy the nice properties of orthogonality for measures; in Example \ref{matrices} we will show that $\ll$ on $S(M_2(\C))$ does not satisfy the ccc-below property and moreover we will construct an analytic mof $\A \sq S(M_2(\C))$. On the other hand, we are still able to use the idea of \cite{strongregularity} to prove that $\sim$ for states is not classifiable by countable structures (this is Proposition \ref{not_classifiable_states}).
    
    There is however an alternative notion of orthogonality (which does not come from a pre-order) and for which it turns out that the argument of \cite{strongregularity}, using product measures, can be generalised to a class of C*-algebras; this is our Theorem \ref{thm_product_states}. Furthermore, we are also able to generalise the result of Kakutani from product measures to product states; this is Proposition \ref{kakutani_states}.
    
    For reader's convenience we start by presenting in detail the definition of absolute continuity for states and the stronger notion of orthogonality. Since we are using descriptive set theoretic methods, we always require $A$ to be separable.
    
    So let $A$ be a separable unital C*-algebra. By the Banach-Alaoglu theorem $(A^*)_1$, the closed unit ball in the dual space of $A$, is compact Polish in the weak*-topology. Moreover, $S(A) \sq (A^*)_1$ is compact convex, so convex compact Polish on its own.
    
    Now let $H$ be the Hilbert space from the universal representation of $A$ in $B(H)$. Recall that we identify $A'' \sq B(H)$, the double commutant of $A$, which is a von Neumann algebra, with the double dual $A^{**}$, and call it the \emph{enveloping von Neumann algebra} (see \cite{brown_ozawa}). Denote $\M := A^{**}$, and recall that we can identify $S(A)$ with $(\M_*)^+_1$, where we identify $\phi \in S(A)$ with its normal extension $\phi^{**} \in S(\M)$, and where $\M_*$ is the predual of $\M$. By definition, $\phi$ being normal means that $\phi$ is ultraweakly continuous on $\M$. We will employ Theorem 7.1.12 of \cite{kadison_ringrose_ii}, saying that a state $\phi$ is normal precisely when it is \emph{completely additive}, i.e., for every orthogonal family of projections $\{p_i\}_{i \in I}$ of $\M$ it holds that $\phi(\sum_{i \in I} p_i) = \sum_{i \in I} \phi(p_i)$. See e.g. \cite{brown_ozawa} or \cite{kadison_ringrose_ii} for more details.
    
    Recall that projections in a von Neumann algebra form a complete complemented lattice with 0 and 1. We use $\vee$ and $\wedge$ to denote \emph{suprema} and \emph{infima} respectfully. So for $\phi \in S(A)$ it holds that
    \[
    \bigvee \{p \in \proj(\M) : \phi(p) = 0\}
    \]
    is a projection in $\M$. Moreover, observe that it also holds that
    \[
    \phi\left(\bigvee \{p \in \proj(\M) : \phi(p) = 0\}\right) = 0.
    \]
    This is because $\vee \{p \in \proj(\M) : \phi(p) = 0\} = \sum_{i \in I} p_i$ for any maximal orthogonal family of projections $\{p_i\}_{i \in I} \sq \M$, satisfying that $\phi(p_i) = 0$ for every $i \in I$. For $\phi \in S(A)$, we define its \emph{support} (sometimes called \emph{carrier}) by
    \[
    \supp \phi := 1 - \bigvee \{p \in \proj(\M) :\phi(p) = 0\},
    \]
    which is again a projection in $\M$. Notice that by definition $\phi(\supp \phi) = 1$ (we think of $\supp \phi$ as the largest projection where $\phi$ is everywhere non-zero). Then for $\phi, \psi \in S(A)$ put
    \[
    \psi \ll \phi \text{ if and only if } \supp \psi \leq \supp \phi,
    \]
    and say that $\psi$ is \emph{absolutely continuous} with respect to $\phi$. Observe that
    \[
    \psi \ll \phi \text{ if and only if } (\forall p \in \proj(\M))\, \phi(p) = 0 \to \psi(p) = 0.
    \]
    Set also $\phi \sim \psi$ if and only if $\phi \ll \psi$ and $\psi \ll \phi$. We continue with the following useful description of absolute continuity for states.
    
    \begin{claim}\label{char_pos}
        For states $\phi, \psi$ it holds that $\psi \ll \phi$ if and only if for every positive $a \in \M$ we have that $\phi(a) = 0$ implies that $\psi(a) = 0$.
    \end{claim}
    
    \begin{proof}
        The direction from right to left is immediate. For the other direction note that for every positive element $a \in \M$ there is a sequence $(\sum_{i = 0}^{k_n} \lambda^n_i p^n_i)_{n \in \N}$, where $\lambda^n_i \in (0, \infty)$ and $p^n_i \in \proj(\M)$ for every $n \in \N$ and $0 \leq i \leq k_n$, so that $\sum_{i = 0}^{k_n} \lambda^n_i p^n_i \leq a$ for every $n \in \N$ and so that $||a - \sum_{i = 0}^{k_n} \lambda^n_i p^n_i|| \to 0$, as $n \to \infty$. Indeed, for $a$ positive (in particular normal), let $\n$ be the Abelian von Neumann subalgebra generated by $a$. Since the statement clearly holds in $\n$ (and the norm on $\n$ is the restriction of the one on $\M$), the same sequence also converges to $a$ from below in $\M$.
        
        Now, if $\phi(a) = 0$, then $\phi(p^n_i) = 0$ for every $n \in \N$ and every $0 \leq i \leq k_n$. So by assumption also $\psi(p^n_i) = 0$ for all $n \in \N$ and all $0 \leq i \leq k_n$. Hence $\psi(a) = 0$, which completes the proof.
    \end{proof}
    
    We proceed with the following result, which tells us that in order to check whether one state is absolutely continuous with respect to the other we do not need to go to the large enveloping von Neumann algebra.
    
    \begin{claim}\label{independent_ambient}
        Let $\phi, \psi \in S(A)$ for a separable unital C*-algebra $A$. Suppose that $\pi: A \to B(K)$ is a faithful representation of $A$ on a Hilbert space $K$, so that $\phi$ and $\psi$ have unique normal extensions to $\n := A'' \sq B(K)$, which we also denote by $\phi$ and $\psi$ respectively. Then $\psi \ll \phi$ if and only if for every $p \in \proj(\n)$ we have that $\phi(p) = 0$ implies that $\psi(p) = 0$.
    \end{claim}
    
    \begin{proof}
        By the universal property of the enveloping von Neumann algebra $\M$, there is a normal *-epimorphism (i.e., a *-homomorphism which is onto) $\alpha: \M \to \n$, which is equal to identity on $A$. Then since every element of $\M$ is a limit of an ultraweakly converging net $\{x_\xi\}$ from $A$, and since $\phi$ and $\psi$ are normal on both $\M$ and $\n$, it holds for every $x \in \M$ that
        \[
        \phi(x) = \phi(\alpha(x)) \quad \text{and} \quad \psi(x) = \psi(\alpha(x)).
        \]
        Suppose first that $\psi \ll \phi$ and that $\phi(p) = 0$ for $p \in \proj(\n)$. Then since $\alpha$ is onto there is some positive $b \in \M$ so that $\alpha(b) = p$. Hence $\phi(b) = 0$ and by Claim \ref{char_pos} also $\psi(b) = 0$, which in turn implies that $\psi(p) = \psi(\alpha(b)) = \psi(b) = 0$.
        
        Conversely, suppose that for every $p \in \proj(\n)$ it holds that $\phi(p)$ implies that $\psi(p) = 0$. Let $q \in  \proj(\M)$ be such that $\phi(q) = 0$. Then $\alpha(q)$ is a projection in $\n$, so $\psi(\alpha(q)) = 0$ and hence also $\psi(q) = \psi(\alpha(q)) = 0$.
    \end{proof}
    
    \begin{remark}\label{supp_to_supp}
        Observe that for any projection $p \in \n$, there is some projection $q\in \M$ for which $\alpha(q) = p$. To see this, take any positive $a \in \M$ for which $\alpha(a) = p$. Then the sequence $(b_n)_n$, defined by
        \[
        b_n := ||a||^{\frac{1}{n}} \left(\frac{a}{||a||}\right)^{\frac{1}{n}},
        \]
        converges (strongly and ultraweakly) to some projection $q \in \M$ (see the proof of Theorem 17.5 from \cite{zhu}), and also satisfies that $\alpha(b_n) = p$ for every $n \in \N$. Hence also $\alpha(q) = p$. So we could have proven the claim without using Claim \ref{char_pos}.
        
        As a consequence of this fact, we actually have that $\alpha(\supp \psi) = \supp \psi$, since clearly
        \[
        \alpha\left(\bigvee \{p \in \proj(\M) : \psi(p) = 0\}\right) \leq  \bigvee \{p \in \proj(\n) : \psi(p) = 0\},
        \]
        but also if $q \in \proj(\M)$ is such that $\alpha(q) = \vee \{p \in \proj(\M) : \psi(p) = 0\}$, then $\psi(q) = 0$.
    \end{remark}
    
    \begin{remark}
        Note that Claim \ref{independent_ambient} holds more generally (with essentially the same proof) for non-degenerate (i.e., unital) representations $\pi: A \to B(K)$, which are not necessarily faithful, and for $\phi, \psi \in S(A)$, for which there are $x, y \in K$, so that for every $a \in A$ it holds that $\phi(a) = \langle \pi(a) x, x \rangle_K$ and $\psi(a) = \langle \pi(a) y, y \rangle_K$.
    \end{remark}
    
     Using (the proof of) Claim \ref{char_pos}, we can replace the requirement in the statement of Claim \ref{independent_ambient} that $p \in \proj(\n)$ with $p \in \n$ being positive.
    
    The reason why we are allowed to call $\ll$ absolute continuity for states is that for commutative C*-algebras it coincides with the classical notion defined for measures. 
    
    \begin{proposition}\label{coincide}
        If $A$ is commutative, $\ll$ defined for states coincides with $\ll$ defined for measures.
    \end{proposition}
    
    \begin{proof}
        Suppose that $A = C(X)$ for a compact Polish space $X$. Let $\mu, \nu \in P(X)$ and let $\phi_\mu, \phi_\nu$ be the corresponding states on $A$ (via the Riesz--Markov--Kakutani representation) and also on $\M = A^{**}$ (via the unique normal extension). Let $\n := L^\infty(X, \frac{1}{2}(\mu + \nu))$ and observe that both $\phi_\mu$ and $\psi_\nu$ admit unique normal extensions to $\n$, again denoted with $\phi_\mu$ and $\psi_\nu$ respectively. Actually, it holds that the extensions to $\n$ are
        \[
        \phi_\mu = \int \,  \diff \mu \quad \text{and} \quad \phi_\nu = \int\,  \diff \nu.
        \]
        Recalling that projections in $\n$ are of the form $\chi_B$ for $B \sq X$ Borel, an application of Claim \ref{independent_ambient} completes the proof.
    \end{proof}
    
    We now move to different notions of orthogonality. As already alluded to above, we say that $\phi, \psi \in S(A)$ are \emph{orthogonal}, which we denote with $\phi \perp \psi$, if there is no $\rho \in S(A)$ for which $\rho \ll \phi$ and $\rho \ll \psi$. In particular, if $\supp \phi \land \supp \psi = 0$, then $\phi \perp \psi$. Hence, it is very easy for two states to be orthogonal. Consequently, $\perp$ for states fails to satisfy the nice properties of $\perp$ for measures.
    
    \begin{example}\label{matrices}
        Let $A := M_2(\C)$. Note that since $A$ is finite dimensional the norm topology coincides with the ultraweak operator topology (and also with the weak/strong/ultrastrong operator topologies). Hence all states on $A$ are normal, so using Claim \ref{independent_ambient}, we can calculate $\ll$ (and $\perp$) by considering supports in the von Neumann algebra $A$. A rank 1 projection $p \in \proj(A)$ (i.e., when $p$ is seen as a projection onto a subspace $E_p$ of $\C^2$, the dimension of $E_p$ is 1), gives rise to a state $\phi_p \in S(A)$, defined for $a \in M_2(\C)$ by
        \[
        \phi_p(a) := \tr(p\, a\, p),
        \]
        where $\tr$ is the usual non-normalised trace on $M_2(\C)$. Note that $\supp \phi_p = p$ and since for two distinct rank 1 projections $p, q$ it holds that $p \land q = 0$, we get that $\phi_p \perp \phi_q$. Hence
        \[
        \A:= \{\phi_p : p \in \proj(A) \text{ has rank } 1\} \sq S(A)
        \]
        is an orthogonal family. Take now any $\psi \in S(A)$. If $\supp \psi = \C^2$, then we get that $\phi \ll \psi$ for all $\phi \in \A$. Since there clearly are such $\psi$ (e.g., take $\psi = \frac{1}{2} (\ev_{1, 1} + \ev_{2, 2})$), the ccc-below property fails for $\ll$. Suppose now that $\supp \psi$ is of rank 1. Then $\psi \sim \phi_{\supp \psi}$. Thus $\A$ is maximal. Moreover, $\A$ is clearly analytic, so we have an analytic mof in $S(A)$.
    \end{example}
    
    Next we introduce another notion of orthogonality. Let $A$ be a separable unital C*-algebra. For $\phi, \psi \in S(A)$, we say that they are \emph{strongly orthogonal}, denoted by $\phi \perpp \psi$ if $\supp \phi \, \supp \psi = 0$. Clearly, $\phi \perpp \psi$ implies $\phi \perp \psi$. Moreover, for $\phi \sim \psi$ and any $\chi \in S(A)$ it holds that $\phi \perpp \chi$ if and only if $\psi \perpp \chi$. We proceed with a useful characterisation of strong orthogonality.
    
    \begin{fact}\label{equivalent_strong}
        For $\phi, \psi \in S(A)$ the following are equivalent
        \begin{enumerate}
            \item $\phi \perpp \psi$;
            \item $\phi(\supp \psi) = 0$;
            \item $(\exists p \in \proj(A^{**}))\, \phi(p) = 0 \land \psi(p) = 1$.
        \end{enumerate}
    \end{fact}
    
    \begin{proof}
        $(1) \Leftrightarrow (2)$ follows from the definition of support and the fact that for projections $p, q$ it holds that $p\, q = 0$ if and only if $p \leq 1-q$.
        
        $(2) \Leftrightarrow (3)$ is clear.
    \end{proof}
    
        With this characterisation, Claim \ref{independent_ambient} and Remark \ref{supp_to_supp}, we obtain the following, which enables us to decide $\perpp$ in smaller representations.
    
    \begin{claim}\label{independent_ambient_perpp}
        Let $\phi, \psi \in S(A)$ for a separable unital C*-algebra $A$. Suppose that $\pi: A \to B(K)$ is a faithful representation of $A$ on a Hilbert space $K$, so that $\phi$ and $\psi$ have unique normal extensions to $\n := A'' \sq B(K)$, which we also denote by $\phi$ and $\psi$ respectively. Then $\phi \perpp \psi$ if and only if $\phi(\supp \psi) = 0$ holds in $\n$.
    \end{claim}
    
    \begin{proof}
        Let $\alpha: \M \to \n$ be as in the proof of Claim \ref{independent_ambient}. By Fact \ref{equivalent_strong}, $\phi \perpp \psi$ if and only if $\phi(\supp \psi) = 0$ holds in $\M$. But by the proof of Claim \ref{independent_ambient} and by Remark \ref{supp_to_supp}, $\phi(\supp \psi) = 0$ is true in $\M$ precisely when it is true in $\n$.
    \end{proof}
    
    \begin{remark}
        As with Claim \ref{independent_ambient}, we actually have that Claim \ref{independent_ambient_perpp} holds more generally for non-degenerate representations $\pi: A \to B(K)$, which are not necessarily faithful, and for $\phi, \psi \in S(A)$, for which there are $x, y \in K$, so that for every $a \in A$ it holds that $\phi(a) = \langle \pi(a) x, x \rangle_K$ and $\psi(a) = \langle \pi(a) y, y \rangle_K$.
    \end{remark}
    
    Notice also that (3) of Fact \ref{equivalent_strong} implies that $\perpp$ is analytic. Next we note that the notion of strong orthogonality extends orthogonality for measures.
    
    \begin{claim}
        Suppose that $A$ is a commutative separable unital C*-algebra and take any $\phi, \psi \in S(A)$. Then $\phi \perpp \psi$ if and only if $\phi \perp \psi$.
    \end{claim}
    
    \begin{proof}
        The forward direction is obvious. So assume that $\phi \perp \psi$ and denote $p := \supp \phi$ and $q := \supp \psi$. Suppose for contradiction that $\phi(q) > 0$. Then define
        $\chi$ by setting
        \[
        \chi(a) := \frac{1}{\phi(q)} \phi(p \, q\, a)
        \]
        for $a \in A$. Clearly, $\chi \in S(A)$ and $\chi \ll \phi, \psi$ which is a contradiction.
    \end{proof}
    
    As opposed to $\perp$ for states, $\perpp$ shares some nice properties with orthogonality for measures. There is a version of the ccc-below property which is true for $\perpp$. A set $\A \sq S(A)$ is a \emph{$\perpp$-antichain} if for every $\psi \neq \chi \in \A$ it holds that $\psi \perpp \chi$.
    
    \begin{claim}\label{ccc_below}
        For any $\phi \in S(A)$ and any $\perpp$-antichain $\A \sq S(A)$ there are only countably many $\psi \in \A$, for which $\neg \phi \perpp \psi$.
    \end{claim}
    
    \begin{proof}
        Suppose for contradiction that there is an uncountable set $\{\psi_i : i \in I\}$ of pairwise strongly orthogonal states, such that for all $i \in I$ it holds that $\neg \phi \perpp \psi_i$. Then $\phi(\supp \psi_i) > 0$ for all $i \in I$. But since $\supp \psi_i$ are pairwise orthogonal and since $\phi$ is completely additive on $A^{**}$, we have that
        \[
        \phi\left(\sum_{i \in I} \supp \psi_i\right) = \sum_{i \in I} \phi(\supp \psi_i).
        \]
        But since $I$ is uncountable, this is a contradiction, as the sum on the right diverges.
    \end{proof}
    
    In \cite{bures}, Bures defined $\rho$ and $d$ for normal states on von Neumann algebras, which generalise the identically denoted notions defined for measures in \cite{kakutani}; see Proposition 2.7 of \cite{bures}. We reintroduce these notions for states on a separable C*-algebra $A$. Let $\phi, \psi \in S(A)$ and put
        \begin{multline*}
            Q(\phi, \psi) := \{(\pi, x, y) : \pi \text{ is a faithful representation of } A \text{ on } H,\\  x, y \in H \text{ respectively induce } \phi, \psi \text{ relative to } \pi \}.
        \end{multline*}
    Here, $x$ \emph{induces} $\phi$ \emph{relative to} $\pi$, means that for every $a \in A$ it holds that $\phi(a) = \langle \pi(a) x, x\rangle_H$. Then define
        \[
        \rho(\phi, \psi) := \sup \{|\langle x, y \rangle| : (\pi, x, y) \in Q(\phi, \psi)\}
        \]
        and
        \[
        d(\phi, \psi) = \inf \{||x - y|| : (\pi, x, y) \in Q(\phi, \psi)\}. 
        \]
        Observe that it holds that $d(\phi, \psi)^2 = 2 (1 - \rho(\phi, \psi))$. The same argument as the one in the proof of Proposition 1.7 from \cite{bures}, shows that $d$, defined for $A$ (not necessarily a von Neumann algebra), is a metric. Let $H_S$ be some separable Hilbert space, for which there is a faithful representation $\pi_S: A \to B(H_S)$. The proof of Proposition 1.6 from \cite{bures} shows that for $\phi, \psi \in S(A)$, we can calculate $\rho(\phi, \psi)$ and $d(\phi, \psi)$ by ranging over the tuples $(\pi, x, y) \in Q(\phi, \psi)$, for which $\pi$ is fixed to be
        \[
        \pi: A \to B(H_\phi \oplus H_\psi \oplus H_S \oplus H_\phi \oplus H_\psi \oplus H_S),
        \]
        defined for $a \in A$ by $\pi(a) := \pi_\phi(a) \oplus \pi_\psi(a) \oplus \pi_S(a) \oplus \pi_\phi(a) \oplus \pi_\psi(a) \oplus \pi_S(a)$, where $\pi_\phi: A \to B(H_\phi)$ and $\pi_\psi: A \to B(H_\psi)$ are the GNS representations. Fixing some countable dense set $D \sq A$, observe that unfolding the definitions of the GNS representations and direct sums of Hilbert spaces, yields that
        \begin{align*}
            \rho(\phi, \psi) = &\sup \{\langle x_0, u_0 \rangle_\phi + \langle y_0, v_0 \rangle_\psi + \langle z_0, w_0 \rangle_{H_S} + \langle x_1, u_1 \rangle_\phi + \langle y_1, v_1 \rangle_\psi + \langle z_1, w_1 \rangle_{H_S} :\\
            &\phantom{\sup a } x_0, x_1, u_0, u_1 \in H_\phi; y_0, y_1, v_0, v_1 \in H_\psi; z_0, z_1, w_0, w_1 \in H_S \text{ and }\\
            &\phantom{\sup a } (\forall d \in D)\, \phi(d) = \langle \pi_\phi(d) x_0, x_0 \rangle_\phi + \langle \pi_\psi(d) y_0, y_0 \rangle_\psi + \langle \pi_S(d) z_0, z_0 \rangle_{H_S}\\
            &\phantom{\sup a (\forall d \in D)\, \phi(d)\, } + \langle \pi_\phi(d) x_1, x_1 \rangle_\phi + \langle \pi_\psi(d) y_1, y_1 \rangle_\psi + \langle \pi_S(d) z_1, z_1 \rangle_{H_S} \, \land  \\
            &\phantom{\sup a (\forall d \in D)\,} \psi(d) = \langle \pi_\phi(d) u_0, u_0 \rangle_\phi + \langle \pi_\psi(d) v_0, v_0 \rangle_\psi + \langle \pi_S(d) w_0, w_0 \rangle_{H_S} \\
            &\phantom{\sup a (\forall d \in D)\, \psi(d)\, }+ \langle \pi_\phi(d) u_1, u_1 \rangle_\phi + \langle \pi_\psi(d) v_1, v_1 \rangle_\psi + \langle \pi_S(d) w_1, w_1 \rangle_{H_S} \}.
        \end{align*}
        In particular, for a fixed $\eps \geq 0$ the sets
        \[
        \{(\phi, \psi) \in S(A)^2 : \rho(\phi, \psi) > \eps\}
        \]
        and
        \[
        \{(\phi, \psi) \in S(A)^2 : d(\phi, \psi) < \eps\}
        \]
        are analytic.
        
        The proof of Lemma 1.2 from Promislow's \cite{promislow} proves that $\phi \perpp \psi$ if and only if $d(\phi, \psi) =\sqrt{2}$. In particular, $\neg \phi \perpp \psi$ is analytic, which since we have already observed that $\phi \perpp \psi$ is analytic, implies that $\perpp$ is Borel.
    
    \subsection{No analytic maximal strongly orthogonal families}\label{product_states}
    
    In this subsection we prove Theorem \ref{thm_product_states}.
    
    In \cite{strongregularity}, Kechris and Sofronidis used the theory of turbulence to prove the following (which is Theorem 3.1 in \cite{strongregularity}). 
    
    \begin{theorem}\label{thm_product_measure}
        For any analytic orthogonal $\A \sq P(2^\N)$, there exists $\alpha \in (0, 1)^\N$ such that $\mu_\alpha \perp \mu$ for every $\mu \in \A$, where
        \[
        \mu_\alpha := \prod_{n \in \N} (\alpha(n) \delta_0 + (1 - \alpha(n)) \delta_1).
        \]
    \end{theorem}
    
    The idea of their proof is to build on Kakutani's \cite{kakutani} and define a continuous map $f: 2^\N \to P(2^\N)$, satisfying that for every $x, y \in 2^\N$ it holds that $x E_I y$ implies that $f(x) \sim f(y)$ and $\neg x E_I y$ implies that $f(x) \perp f(y)$, where
    \[
    x E_I y \quad \text{if and only if} \quad \sum_{n \in x \Delta y} \frac{1}{n} < \infty.
    \]
    (We use the notation $E_I$ because it is used in \cite{strongregularity}; other more common notations are either $I_2$ or $E_2$.) Recall that a Borel equivalence relation $E$ on a Polish space $Y$ is \emph{generically $S_\infty$-ergodic} if every $E$-class is meagre and for any standard Borel space $Z$, equipped with a Borel action of $S_\infty$, and any Baire measurable $f: Y \to Z$, with the property that $x E y$ implies that $(\exists g \in S_\infty)\, g \cdot f(x) = f(y)$, there is an $E$-invariant comeagre set $C \sq Y$, such that $f$ maps $C$ to a single class in $Z$.
    
    If $E$ is generically $S_\infty$-ergodic, then $E$ is not classifiable by countable structures. Here a relation $E$ on a standard Borel space $X$ is said to be \emph{classifiable by countable structures} if there is a countable language $L$ and a Borel map $f: X \to X_L$ (where $X_L$ is the space of countable structures for $L$), so that for all $x, y \in X$ it holds that $x E y$ if and only if $f(x) \cong f(y)$. See \cite{hjorth} for more about classification by countable structures, generic ergodicity and turbulence.
    
    The relation $E_I$ defined above is generically $S_\infty$-ergodic (see \cite{hjorth} and \cite{strongregularity}) and so with the above reduction of $E_I$ to $\sim$, Kechris and Sofronidis establish that $\sim$ is not classifiable by countable structures. Then they prove the following lemma (see Lemma 3.3 in \cite{strongregularity}), which gives Theorem \ref{thm_product_measure} (since $\ll$ for measures has the ccc-below property).
    
    \begin{lemma}\label{S_ergodic}
        Let $\preccurlyeq$ be an analytic partial pre-ordering on a Polish space $X$ which satisfies the ccc-below property and assume that there exists a generically  $S_\infty$-ergodic equivalence relation $E$ on a Polish space $Y$ and a Borel measurable function $f: Y \to X$ with the properties that $z E y \implies f (z) \sim f (y)$ and $\neg z E y \implies f (z) \perp f (y)$, whenever $z, y$ are in $Y$. Then, given any orthogonal analytic subset $\A$ of X, there exists $y \in Y$ such that $f(y) \perp x$ for every $x \in \A$.
    \end{lemma}
    
    We will prove that the same idea can be generalised to a class of C*-algebras. 
    Recall that for a C*-algebra $A$ and its subalgebra $B$, a linear map $E: A \to B$ is a \emph{conditional expectation}, when $E$ is a contractive completely positive projection, such that for every $a \in A$ and $b, b' \in B$ it holds that $E(b a b') = b E(a) b'$. By Tomiyama's theorem (see Theorem 1.5.10 in \cite{brown_ozawa}), a projection $E: A \to B$ is a conditional expectation precisely when it is contractive. Notice that conditional expectations are closed under composition.
    
    Suppose now that $B \sq A$ are unital C*-algebras (with possibly different units) and that $E: A \to B$ is a conditional expectation (note that $E(1_A) = 1_B E(1_A) 1_B = E(1_B) = 1_B$). For a state $\phi \in S(B)$, there is an extension $\Tilde{\phi} \in S(A)$, defined by $\tilde{\phi}(a) = \phi(E(a))$ for $a \in A$. Clearly the map $\Tilde{(\cdot)}: S(B) \to S(A)$, defined by $\phi \mapsto \Tilde{\phi}$ is continuous. Note also that $E^{**}: A^{**} \to B^{**}$ is again a conditional expectation, extending $E$. We next list some nice properties of $\tilde{(\cdot)}$.
    
    \begin{claim}\label{preserve_continuity}
        For $\phi, \psi \in S(B)$ it holds that $\psi \ll \phi$ if and only if $\Tilde{\psi} \ll \Tilde{\phi}$.
    \end{claim}
    
    \begin{proof}
        Suppose that $\psi \ll \phi$ and that $\Tilde{\phi}(a) = 0$ for a positive $a \in A^{**}$. This means that $\phi(E^{**}(a)) = 0$, and consequently $\psi(E^{**} (a)) = 0$. Hence by definition, $\Tilde{\psi}(a) = 0$.
        
        The other direction is obvious.
    \end{proof}
    
    \begin{claim}\label{preserve_orthogonality}
        For $\phi, \psi \in S(B)$ it holds that $\phi \perpp \psi$ implies $\tilde{\phi} \perpp \tilde{\psi}$.
    \end{claim}
    
    \begin{proof}
        Suppose that $\phi \perpp \psi$. We will show that $\tilde{\phi}(\supp \tilde{\psi}) = 0$, which implies that $\tilde{\phi} \perpp \tilde{\psi}$ by Fact \ref{equivalent_strong}. Note that
        \[
        \supp \tilde{\psi} = 1_A - \bigvee \{p \in \proj(A^{**}) : \psi(E^{**}(p)) = 0\}.
        \]
        Since
        \[
        \bigvee \{q \in \proj(B^{**}) : \psi(q) = 0\} \leq \bigvee \{p \in \proj(A^{**}) : \psi(E^{**}(p)) = 0\},
        \]
        and since $E^{**}$ is monotone, it holds that
        \[
        E^{**}(\supp \tilde{\psi}) \leq 1_B - \bigvee \{q \in \proj(B^{**}) : \psi(q) = 0\} = \supp \psi.
        \]
        But then
        \[
        \tilde{\phi}(\supp \tilde{\psi}) = \phi(E^{**}(\supp \tilde{\psi})) \leq \phi(\supp \psi) = 0
        \]
        and hence $\tilde{\phi}(\supp \tilde{\psi}) = 0$, which completes the proof.
    \end{proof}
    
    With this we are ready to prove Theorem \ref{thm_product_states}.
    
    \begingroup
    \def\thetheorem{\ref{thm_product_states}}
    \begin{theorem}
        Suppose $A$ is a separable unital C*-algebra, which contains a copy of $C(2^\N)$ as a subalgebra and for which there is a conditional expectation $E: A \to C(2^\N)$. Then for every strongly orthogonal $\A \sq S(A)$ there is $\alpha \in (0, 1)^\N$ so that $\Tilde{\mu}_\alpha \perpp \psi$ for every $\psi \in \A$, where $\Tilde{\mu}_\alpha$ is the extension of the state, corresponding to the product measure
        \[
        \prod_{n \in \N} (\alpha(n) \delta_0 + (1-\alpha(n)) \delta_1),
        \]
        from $C(2^\N)$ to $A$.
    \end{theorem}
    \addtocounter{theorem}{-1}
    \endgroup
    
    \begin{proof}
        As in \cite{strongregularity}, define $\alpha: 2^\N \to [\frac{1}{4}, \frac{3}{4}]^\N$ by
        \[
        \alpha(x)(n) := \begin{cases}
        \frac{1}{4}\left(1 + \frac{1}{\sqrt{n + 1}}\right) &\text{if } n \in x\\
        \frac{1}{4} &\text{if } n \in \N \setminus x
        \end{cases}
        \]
        for all $x \in 2^\N$, where we identify $2^\N$ with $\pp(\N)$, the powerset of $\N$. Now let $f: 2^\N \to S(A)$ be defined by $f(x) = \Tilde{\mu}_{\alpha(x)}$ for $x \in 2^\N$. Since the maps $\alpha$, $\Tilde{(\cdot)}$ and the map $[\frac{1}{4}, \frac{3}{4}]^\N \to S(2^\N)$, defined by $h \mapsto \mu_h$ are all continuous, so is $f$.
        
        In \cite{strongregularity}, it is established that for $x, y\in 2^\N$ it holds that
        \[
        \sum_{n \in \N} (\alpha(x)(n) - \alpha(y)(n))^2 = \sum_{n \in x \Delta y} \frac{1}{16 (n + 1)},
        \]
        so by Corollary 1 from \cite{kakutani} and by Claims \ref{preserve_continuity} and \ref{preserve_orthogonality} we have for every $x, y\in 2^\N$ that
        \[
        x E_I y \implies f(x) \sim f(y) \quad \text{and} \quad \neg x E_I y \implies f(x) \perpp f(y).
        \]
        Lemma \ref{S_ergodic} still holds (with the same proof) with ccc-below replaced with the property from Claim \ref{ccc_below} and with $\preccurlyeq$ being analytic replaced with $\perpp$ co-analytic (in our case it is even Borel). Thus the proof is complete.
    \end{proof}
    
    The function $f$ from the above proof is a continuous reduction of $E_I$ to $\sim$ for states, and hence we have the following consequence.
    
    \begin{corollary}\label{not_classifiable_states}
        Suppose $A$ is a separable unital C*-algebra, which contains a copy of $C(2^\N)$ as a subalgebra and for which there is a conditional expectation $E: A \to C(2^\N)$. Then $\sim$ on $S(A)$ is not classifiable by countable structures.
    \end{corollary}
    
    The following examples of nice C*-algebras satisfying assumptions of Theorem \ref{thm_product_states} (and hence also of Corollary \ref{not_classifiable_states}) were suggested to the author by Magdalena Musat and Mikael R{\o}rdam.
    
    Fix any $A$ satisfying assumptions of Theorem \ref{thm_product_states}. Let $B$ be any separable unital C*-algebra and pick some $\phi \in S(A \otimes B)$. Then there is a conditional expectation $E: A \otimes B \to A$, induced by
    \[
    E(a \otimes b) = \phi(b)\, a.
    \]
    Hence $A \otimes B$ also satisfies assumptions of Theorem \ref{thm_product_states}.
    
    Recalling Proposition 4.1.9 from \cite{brown_ozawa}, we get that for any countable discrete group $\Gamma$ and any homomorphism $\alpha: \Gamma \to \aut(A)$, the reduced crossed product $A \rtimes_{\alpha, r} \Gamma$ also satisfies the assumptions of Theorem \ref{thm_product_states}. Note that if $A$ is simple and the action is \emph{outer} (i.e., $\Gamma$ acts by outer automorphisms), then by \cite{kishimoto} $A \rtimes_{\alpha, r} \Gamma$ is simple. 
    
    Consider the following diagram
    \[
    \begin{tikzcd}
    \C \ar["{\phi_0}", r] \ar["E_0", d] & M_2(\C) \ar["{\phi_1}", r] \ar["E_1", d] & M_4(\C) \ar["{\phi_2}", r] \ar["E_2", d] & M_8(\C) \ar["{\phi_3}", r] \ar["E_3", d]& \cdots \ar[r]& M_{2^\infty} \ar["E", d]\\
    \C \ar["{\psi_0}", r] & \C^2 \ar["{\psi_1}", r] & \C^4 \ar["{\psi_2}", r] & \C^8 \ar["{\psi_3}", r] & \cdots \ar[r] & C(2^\N),
    \end{tikzcd}
    \]
    where $\phi_n: M_{2^n}(\C) \to M_{2^{n+1}}(\C)$ is defined as $\phi_n(X) := \id_2 \otimes X$, $\psi_n: \C^{2^n} \to \C^{2^{n+1}}$ as $\psi_n(a_1, a_2,  \ldots, a_{2^n}) := (a_1, a_1, a_2, a_2, \ldots, a_{2^n}, a_{2^n})$ and $E_n: M_{2^n}(\C) \to \C^{2^n}$ as $E([a_{i, j}]) = (a_{1, 1}, a_{2, 2}, \ldots, a_{2^n, 2^n})$. On the right we have the inductive limits of the respective sequences: the CAR algebra $M_{2^\infty}$ and $C(2^\N)$. Note that all $E_n$ are conditional expectations and that the diagram commutes, which gives a conditional expectation $E: M_{2^\infty} \to C(2^\N)$. Since  $M_{2^\infty}$ is also separable, it satisfies the assumptions of Theorem \ref{thm_product_states}. 
    
    Since one can view the Cuntz algebra $\mathcal{O}_2$ as a crossed product of $M_{2^\infty}$ with integers (see \cite{cuntz}, \cite{fermion} and \cite{cuntz_crossed}), $\mathcal{O}_2$ also satisfies the assumptions of Theorem \ref{thm_product_states}.
   
   \subsection{Kakutani's theorem for states}
   
   In this subsection we use results from the previous subsection to prove Proposition \ref{kakutani_states}.
   
    Recall that we may view the CAR algebra $M_{2^\infty}$ as $\otimes_{n = 0}^{\infty} M_2(\C)$. If $(\phi_n)_{n \in \N}$ is a sequence of states on $M_2(\C)$, then $\otimes_{n = 0}^\infty \phi_n$ denotes the unique state (called \emph{product state}) with the property that for every sequence $(a_n)_{n \in \N}$, where for all but finitely many $n$ it holds that $a_n = 1_{M_2(\C)}$, we have that
    \[
    \left(\bigotimes_{n = 0}^\infty \phi_n\right) \left(\bigotimes_{n = 0}^\infty a_n \right) = \prod_{n \in \N} \phi_n(a_n).
    \]
    Observe that $\ev_{1, 1}$ and $\ev_{2, 2}$ are states on $M_2(\C)$. These will be our non-commutative analogues of the Dirac measures $\delta_0$ and $\delta_1$ on $2 = \{0, 1 \}$, used in \cite{kakutani} and \cite{strongregularity}.
    
    In \cite{bures}, Bures extended Kakutani's theorem to semi-finite von Neumann algebras. This was improved by Promislow (see \cite{promislow}) to general von Neumann algebras. However, their statements do not mention absolute continuity nor (strong) orthogonality for states, which on the other hand are central to Kakutani's statement. Using Claims \ref{preserve_continuity} and \ref{preserve_orthogonality} and Corollary 1 of \cite{kakutani}, we provide an extension of Kakutani's result about absolute continuity and orthogonality to the special case for product states on the CAR algebra.
    
    \begingroup
    \def\thetheorem{\ref{kakutani_states}}
    \begin{proposition}
    Suppose that $(\alpha_n)_{n \in \N}, (\beta_n)_{n \in \N} \in [\frac{1}{4}, \frac{3}{4}]^\N$ and let
    \[\phi_n := \alpha_n \ev_{1, 1} + (1-\alpha_n) \ev_{2, 2} \quad \text{and} \quad  \psi_n := \beta_n \ev_{1, 1} + (1 -\beta_n) \ev_{2, 2}
    \]
    be states on $M_2(\C)$. Let also $\phi := \otimes_{n = 0}^\infty \phi_n$ and $\psi := \otimes_{n = 0}^\infty \psi_n$ be the product states on $M_{2^\infty}$. Then in $S(M_{2^\infty})$, either $\phi \sim \psi$ or $\phi \perpp \psi$ according to whether
    \[
    \sum_{n \in \N} (\alpha_n - \beta_n)^2
    \]
    converges or diverges respectively.
    \end{proposition}
    \addtocounter{theorem}{-1}
    \endgroup
    
    \begin{proof}
        Let
        \[
        \mu := \prod_{n \in \N} \alpha_n \delta_0 + (1-\alpha_n) \delta_1 \quad \text{and} \quad \nu := \prod_{n \in \N} \beta_n \delta_0 + (1-\beta_n) \delta_1
        \]
        be product measures and note that $\phi = \tilde{\mu}$ and $\psi = \tilde{\nu}$ (using the conditional expectation $E$ defined in the end of the previous subsection).
        
        Now by Corollary 1 of \cite{kakutani}, we get that either $\mu \sim \nu$ or $\mu \perp \nu$ (which is the same as $\mu \perpp \nu$) according to whether
        \[
        \sum_{n \in \N} (\alpha_n - \beta_n)^2
        \]
        converges or diverges respectively. But then by Claims \ref{preserve_continuity} and \ref{preserve_orthogonality}, we get that $\phi \sim \psi$ if and only if $\mu \sim \nu$ and $\phi \perpp \psi$ if and only if $\mu \perp \nu$, which completes the proof.
    \end{proof}
    
    \begin{remark}
        The main ingredient of \cite{kakutani} are properties of $\rho$ and $d$, defined for measures. Since for every $a \in M_2(\C)$ it holds that
        \[
        \phi_n(a) = \tr\left(\begin{pmatrix}
        \sqrt{\alpha_n} & 0\\
        0 & \sqrt{1 - \alpha_n}
        \end{pmatrix} a \begin{pmatrix}
        \sqrt{\alpha_n} & 0\\
        0 & \sqrt{1 - \alpha_n}
        \end{pmatrix}
        \right)
        \]
        and
        \[
        \psi_n(a) = \tr\left(\begin{pmatrix}
        \sqrt{\beta_n} & 0\\
        0 & \sqrt{1 - \beta_n}
        \end{pmatrix} a \begin{pmatrix}
        \sqrt{\beta_n} & 0\\
        0 & \sqrt{1 - \beta_n}
        \end{pmatrix}
        \right),
        \]
        Proposition 2.3 from \cite{bures} implies that
        \[
        \rho(\phi_n, \psi_n) = \tr\left(\begin{pmatrix}
        \sqrt{\alpha_n} & 0\\
        0 & \sqrt{1 - \alpha_n}
        \end{pmatrix} \begin{pmatrix}
        \sqrt{\beta_n} & 0\\
        0 & \sqrt{1 - \beta_n}
        \end{pmatrix}  \right) = \sqrt{\alpha_n \beta_n } + \sqrt{(1- \alpha_n) (1-\beta_n)},
        \]
        which is equal to $\rho(\mu_n, \nu_n)$ (see paragraph above Corollary 1 in \cite{kakutani}). So we actually have that either $\phi \sim \psi$ or $\phi \perpp \psi$ according to whether $\prod_{n \in \N} \rho(\phi_n, \psi_n)$ is positive or equal to 0. Thus in absence of Claims \ref{preserve_continuity} and \ref{preserve_orthogonality}, one could prove Proposition \ref{kakutani_states} by combining the proofs of \cite{kakutani} and \cite{bures} (and \cite{promislow}).
    \end{remark}
    
    \section{Conclusion and open problems}\label{open_problems}
    
    We conclude the paper with discussions about related topics and open questions.
    
    \subsection{Abstract theorem}
    
    The proof of Theorem \ref{thm_measures_comeagre} can be used to prove a more general fact.
    
    \begin{theorem}\label{thm_abstract}
        Suppose that there is a semi-normed vector space $(E, || \cdot ||)$, which has a convex subset $X$, contained in the closed unit ball of $E$, so that for any $x, y \in X$ it holds that $|| x - y|| \leq 1$. Moreover, $X$ carries a Polish topology $\tau$, which has a basis consisting of convex sets, so that for any $\eps \in (0, \infty)$ the set $\{(x, y) \in X \times X : ||x - y|| < \eps \}$ is analytic with respect to $\tau$. A subset $A \sq X$ is called an \emph{antichain} if for any two $y \neq z \in A$ we have that $||y - z || = 1$. Suppose finally that the following properties are satisfied:
        \begin{enumerate}
            \item for every $x \in X$, the set $\{y \in X : ||x - y|| < 1\}$ does not contain an uncountable antichain;
            \item
            for every $\eps \in (0, 1)$, $x \in X$ and any antichain $A \sq X$ the set $\{y \in A : ||x - y|| < \eps \}$ is finite;
            \item for every $x \in X$ the set $x^\perp := \{ y \in X : || x - y || = 1\}$ is comeagre in $(X, \tau)$.
        \end{enumerate}
        Then for any analytic antichain $\A$, it holds that
        \[
        \A^\perp := \{x \in X : (\forall y \in \A) \, ||x - y|| = 1\}
        \]
        is comeagre.
    \end{theorem}
    
    \begin{proof}[Sketch of proof.]
        Suppose for contradiction that $\A$ is an analytic antichain, for which $\A^\perp$ is not comeagre. Then there is a non-empty convex open $O \sq X$ in which $\A^\perp$ is meagre. Let $Z \sq O \setminus \A^\perp$ be a dense (in $O$) $G_\delta$ set.
        
        For $k \in \N$ set $E_k$ to be the space of $k$-element subsets of $\A$ and for fixed $\eps \in (0, 1)$ and $\tau \in (0, \eps)$ define
        \[
        H_{k, \eps} := \{x \in Z : (\exists F \in E_k)\, (\forall y \in F) ||x - y|| < \eps\}
        \]
        and
        \[
        U^\tau_{k, \eps} := H_{k, \eps - \tau} \setminus H_{k + 1, \eps},
        \]
        which have the property of Baire. By the same reasoning as for measures, we get that
        \[
        Z = \bigcup_{k \geq 1} \bigcup_{n > 1} \bigcup_{m > n} U^{1/m}_{k, 1/n}
        \]
        and so for some $k, \eps := 1/n$ and $\tau := 1/m$ it holds that $U^\tau_{k, \eps}$ is comeagre in a nonempty open convex set $V \sq O$. Then by an application of Kuratowski-Ulam theorem (see the proof of Claim \ref{Kuratowski-Ulam}), there is $x \in U^\tau_{k, \eps} \cap V$ and a comeagre $C \sq U^\tau_{k, \eps} \cap V$ so that for every $y \in C$ the set
        \[
        M_y := \{t \in [0, 1] : t x + (1-t) y \in U^\tau_{k, \eps}\}
        \]
        is comeagre in $[0, 1]$. Defining for $y \in U^\tau_{k, \eps}$ the set
        \[
        N_y := \{z \in \A : ||y - z|| < \eps - \tau\},
        \]
        we observe as in the proof of Claim \ref{constant} that for every $z \in C$ it holds that $N_x = N_z$. Then let $y_0, \ldots, y_{k-1}\in \A$ be such that $N_x = \{y_0, \ldots, y_{k-1}\}$. Since for every $y \in X$ it holds that $y^\perp$ is comeagre, we get a contradiction, since
        \[
        B:= \bigcap_{i = 0}^{k-1} y_i^\perp
        \]
        and $C$ are both comeagre in $V$.
    \end{proof}
    
    \begin{remark}\label{automatic_3}
        Note that if $X \times X \to \R_{\geq 0}$, defined by $(x, y) \mapsto || x-y||$ is lower-semicontinuous with respect to $\tau$ (so that in particular $\perp$ is $G_\delta$) and if the space of extreme points of $X$ (which is Polish by Proposition 2.1 from \cite{poulsen_simplex}) is perfect and an antichain, then item (3) follows by an argument similar to the proof of Proposition 4.1 from \cite{strongregularity}.
    \end{remark}
    
    \begin{remark}
        Assuming $\MA + \neg \CH$, $\PD$ or $\AD$ we can replace ``analytic'' from the statement of Theorem \ref{thm_abstract} with ``$\mathbf{\Sigma^1_2}$'', ``projective'' or ``any'' respectively.
    \end{remark}
    
    For Theorem \ref{thm_abstract} to have any value, other examples than measures are needed.
    
    \begin{question}
        Are there other natural examples beside Borel probability measures which satisfy assumptions of Theorem \ref{thm_abstract}?
    \end{question}
    
    In the next subsection we discuss a possible candidate.
    
    \subsection{Comeagreness of witnesses to non-maximality for states}
    
    Let $A$ be a separable unital C*-algebra. As we have seen above, there is a metric $d$ on $S(A)$, defined by Bures in \cite{bures}, which satisfies that for $\phi, \psi \in S(A)$ it holds that $\phi \perpp \psi$ if and only if $d(\phi, \psi) = \sqrt{2}$.
    
    This gives us hope that Theorem \ref{thm_abstract} might apply to $\frac{1}{\sqrt{2}} d$. However,  we do not know whether items (2) and (3) of Theorem \ref{thm_abstract} are satisfied.
    Even if this approach fails, it might still be the case that the following has a positive answer. 
    
    \begin{question}
        Is it the case for a separable unital C*-algebra $A$ (or for a class of C*-algebras satisfying some additional properties), that for any analytic strongly orthogonal family $\A \sq S(A)$, the set \[\A^{\perpp} := \{\psi \in S(A) : (\forall \phi \in \A)\, \psi \perpp \phi\}\] is comeagre?
    \end{question}
    
    \subsection{Measure on the space of measures}

    For this subsection we work on the Cantor space $2^\N$. Observe that $P(2^\N)$ is homeomorphic to
    \[
    p(2^\N) := \{f \in [0, 1]^{2^{<\N}} : f(\emptyset) = 1 \land (\forall s \in 2^{<\N})\, f(s) = f(s^\smallfrown 0) + f(s^\smallfrown 1) \}.
    \]
    Actually, there is even an isometric bijection between the two spaces when one defines natural metrics on both spaces, which generate the respective Polish topologies, see \cite{fischer_tornquist}.
    
    Furthermore, there is a surjective continuous map $\Phi: [0, 1]^{2^{<\N}} \to p(2^\N)$, defined recursively by
    \begin{align*}
        \Phi(f)(\emptyset) &:= 1\\
        \Phi(f)(s^\smallfrown 0) &:= \Phi(f)(s)\cdot f(s)\\
        \Phi(f)(s^\smallfrown 1) &:= \Phi(f)(s)\cdot (1 - f(s))
    \end{align*}
    for $f \in [0, 1]^{2^{<\N}} $ and $s \in 2^{<\N}$. Let $\lambda$ denote the Lebesgue measure on $[0, 1]$. Then $\prod_{s \in 2^{< \N}} \lambda$ is a Borel probability measure on $[0, 1]^{2^{<\N}}$ and $\Phi$ is injective on a set of measure 1. We denote the pushforward of this measure to $P(2^\N)$ (via the identifications above) with $\Lambda$. Given Theorem \ref{thm_measures_comeagre}, it is natural to ask the following question.
    
    \begin{question}\label{q_measure_measures}
        Suppose that $\A \sq P(2^\N)$ is an analytic orthogonal family. Does it hold that $\Lambda(\A^\perp) > 0$?
    \end{question}
    
    One should not get one's hopes too high and wonder whether it could be that $\Lambda(\A^\perp) = 1$, as this turns out to be false.
    
    \begin{claim}
        For any $\mu \in P(2^\N)$, $\Lambda(\mu^\perp) < 1$.
    \end{claim}
    
    \begin{proof}[Sketch of proof.]
        Fix any $\mu \in P(2^\N)$ and suppose for contradiction that $\Lambda(\mu^\perp) = 1$. Fix also some arbitrary small $\eps \in (0, 1/2)$. For $s \in 2^{<\N}$, let $U_s := \{x \in 2^\N : s \sq x\}$. By repeated use of Fubini's theorem and the fact that for $B \sq [0, 1]$, $\lambda(B) = 1$ implies that $B$ is dense in $[0, 1]$, we get that there is some $\nu \in \mu^\perp$ so that for all $s \in 2^{< \N}$ it holds that
        \[
        \nu(U_s) \in ((1 - \eps)\, \mu(U_s), (1 + \eps)\, \mu(U_s)).
        \]
        Since open sets are disjoint unions of basic open sets, the same holds for all open $U \sq 2^\N$. But then $\nu$ and $\mu$ are not orthogonal, which is a contradiction.
    \end{proof}
    
    \subsection{Definable maximal orthogonal families in forcing extensions}

    The original motivation for trying to find a short and simple proof of Theorem \ref{thm_measures} was that maybe a new proof would help us answer the following open question (which is a reformulation of Open problem 1) from \cite{mrl}).
    
    \begin{question}
        Are there any $\Pi^1_1$ maximal orthogonal families $\A \sq P(2^\N)$ in Laver extensions?
    \end{question}
    
    The hope was also that a new proof of Theorem \ref{thm_measures}, would shed some light onto why some arboreal forcing notions (see \cite{Lwe1998UniformUA} for the definition and results on arboreal forcing) admit $\Pi^1_1$ mofs in their forcing extensions (Sacks and Miller forcing) and some do not (Mathias forcing); see \cite{mrl}.
    
    \subsection{Nice and bad subsets}\label{nice_bad}

    Let $X$ be a Polish space and $A$ a separable unital C*-algebra. Call an analytic subset $Y \sq P(X)$ (respectively $Y \sq S(A)$) \emph{nice}, if for every analytic pairwise orthogonal (respectively strongly orthogonal) family $\A \sq Y$, there is $\mu \in Y \cap \A^\perp$. Otherwise, call $Y$ \emph{bad}.
    
    For example Theorems \ref{thm_measures} and \ref{thm_product_states} imply that $P(X)$ and $S(A)$ (where $A$ satisfies assumptions of Theorem \ref{thm_product_states}) are nice. Moreover, when $X$ is compact perfect Theorem \ref{thm_measures_comeagre} implies that all non-meagre analytic $Y \sq P(X)$ are nice and Theorems \ref{thm_product_measure} and \ref{thm_product_states} imply that the sets of product measures and product states are nice.
    
    On the other hand $\extr P(X)$ is clearly bad, since $\extr P(X)$ is an orthogonal family (by Proposition 4.1 of \cite{strongregularity}).
    
    \begin{question}\label{q_nice_bad}
        Are there other natural examples of nice/bad sets?
    \end{question}
    
    \bibliographystyle{alpha}
    \bibliography{bibliography}
    
\end{document}